\documentclass{article}

\title{\LARGE \textbf{A Size Bound for Hamilton Cycles}}
\author{Zh.G. Nikoghosyan}

\begin{document}

\maketitle

\begin{abstract}
Every graph of size $q$ (the number of edges) and minimum degree $\delta$   is hamiltonian if $q\le\delta^2+\delta-1$. The result is sharp. \\

\end{abstract}

\section{Introduction}

The earliest two sufficient conditions for a graph to be hamiltonian are based on three simplest graph invariants, namely order $n$, size $q$ and minimum degree $\delta$, in forms of  simple algebraic relations between $n, \delta$ and $n, q$, respectively.\\

\noindent\textbf{Theorem A} (Dirac, 1952)  \cite{[3]}. Every graph with $\delta\ge\frac{n}{2}$ is hamiltonian.\\

\noindent\textbf{Theorem B} (Erd\"{o}s and Gallai, 1959) \cite{[4]}. Every graph with $q\ge\frac{n^2-3n+5}{2}$ is hamiltonian.\\

In this paper we present an analogous simple relation between $\delta$ and q.\\

\noindent\textbf{Theorem 1.} Every graph with $q\le\delta^2+\delta-1$ is hamiltonian.\\

The bound $\delta^2+\delta-1$ in Theorem 1 can not be relaxed to $q\le\delta^2+\delta$ since the graph $K_1+2K_\delta$ consisting of two copies of $K_{\delta+1}$ and having exactly one vertex in common, has $\delta^2+\delta$ edges and is not hamiltonian.

\section{Notations and preliminaries}

Only finite undirected graphs without loops or multiple edges are considered. We reserve $n$, $q$, $\delta$ and $\kappa$  to denote the number of vertices (order), the number of edges (size), the minimum degree and connectivity of a graph. A good reference for any undefined terms is \cite{[1]}. 

The set of vertices of a graph $G$ is denoted by $V(G)$ and the set of edges by $E(G)$. The neighborhood of a vertex $x\in V(G)$ will be denoted by $N(x)$. Set $d(x)=|N(x)|$.  For a subgraph $H$ of $G$ we use $d_H(x)$ short for $|N(x)\cap V(H)|$. Further, we will use $G_m$ to denote an arbitrary graph on $m$ vertices. 

A simple cycle (or just a cycle) $C$ of length $t$ is a sequence $v_1v_2...v_tv_1$ of distinct vertices $v_1,...,v_t$ with $v_iv_{i+1}\in E(G)$ for each $i\in \{1,...,t\}$, where $v_{t+1}=v_1$. When $t=2$, the cycle $C=v_1v_2v_1$ on two vertices $v_1, v_2$ coincides with the edge $v_1v_2$, and when $t=1$, the cycle $C=v_1$ coincides with the vertex $v_1$. So, all vertices and edges in a graph can be considered as cycles of lengths 1 and 2, respectively. A graph $G$ is hamiltonian if $G$ contains a Hamilton cycle, i.e. a cycle of length $n$.  

In order to prove the main result we need Theorem A \cite{[3]} and the following theorem due to Chie \cite{[2]}.\\

\noindent\textbf{Theorem C} (Chie, 1980) \cite{[2]}. Let $G$ be a 2-connected graph. If $d(u)+d(v)\ge n-1$ for each pair of nonadjacent vertices $u,v$ then either $G$ is hamiltonian or  $G=G_{(n-1)/2}+\overline{K}_{(n+1)/2}$.

\section{Proof of Theorem 1}

Assume the converse, that is $G$ satisfies the condition 
$$
q\le \delta^2+\delta-1   \eqno{(1)}
$$
and is not hamiltonian. Since 
$$
q=\frac{1}{2}\sum_{u\in V(G)}d(u)\ge \frac{\delta n}{2},
$$
we have $\delta n/2 \le \delta^2+\delta-1$ which is equivalent to 
$$
\delta\ge \frac{n-1}{2}-\frac{1}{2}+\frac{1}{\delta}.
$$
If $n$ is even, i.e. $n=2t$ for some integer $t$, then 
$$
\delta\ge \frac{2t-1}{2}-\frac{1}{2}+\frac{1}{\delta}=t-1+\frac{1}{\delta},
$$
implying that $\delta\ge t=n/2$. By Theorem A, $G$ is hamiltonian, a contradiction. Let $n$ is odd, i.e. $n=2t+1$ for some integer $t$. Then $\delta\ge t-1/2+1/\delta$ implying that $\delta\ge t\ge (n-1)/2$. Since $G$ is hamiltonian when $\delta > (n-1)/2$ (by Theorem A), we can assume that $\delta=(n-1)/2$. \\

Case 1. $\kappa\ge2$.

By Theorem C, $G=G_{(n-1)/2}+\overline{K}_{(n+1)/2}=G_\delta+\overline{K}_{\delta+1}$. 
Clearly $|E(G)|\ge\delta(\delta+1)$, contradicting (1).\\

Case 2. $\kappa\le1$.

It follows that $G$ has a cut vertex $v$. Let  $H_1$ and $H_2$ be any two connected components of $G\backslash v$. Denote by $H^*_i$ the subgraph induced by $V(H_i)\cup \{v\}$ $(i=1,2)$. Clearly $|V(H^*_i)|\ge \delta+1$ $(i=1,2)$ and $n\ge2\delta+1$. \\

Case 2.1. $xv\not\in E(G)$ for some $x\in V(H_i)$ and $i\in\{1,2\}$.

It follows that $|V(H^*_i)|\ge \delta+2$ $(i=1,2)$ and $n\ge2\delta+2$. Hence
$$
q=\frac{1}{2}\sum_{u\in V(G)}d(u)\ge \frac{\delta n}{2}\ge \delta(\delta+1),
$$
contradicting (1).\\

Case 2.2. $xv\in E(G)$ for each $x\in V(H_i)$ and $i\in \{1,2\}$.

It follows that $d_{H_i}(v)\ge |V(H_i)|\ge \delta$ for $i=1,2$ implying that $d(v)\ge d_{H_1}(v)+d_{H_2}(v)\ge 2\delta$. Hence
$$
q= \frac{1}{2}\sum_{u\in V(G)}d(u)= \frac{1}{2}\left (\sum_{u\in V(G)\backslash v}d(u)+d(v)\right )              
$$
$$
\ge\frac{1}{2}\sum_{u\in V(G\backslash v)}d(u)+\delta\ge \frac{1}{2}(n-1)\delta+\delta\ge\delta^2+\delta, 
$$
contradicting (1).  Theorem 1 is proved.

\noindent Institute for Informatics and Automation Problems\\ National Academy of Sciences\\
P. Sevak 1, Yerevan 0014, Armenia\\  E-mail: zhora@ipia.sci.am

\end{document}